\documentclass[12pt,twoside]{article}
\usepackage{amsmath,amssymb,amsthm,amscd,graphicx,color,accents,algorithm}
\usepackage{psfrag,epsf}
\usepackage{enumerate}

\usepackage{enumitem}

\usepackage{url} 

\usepackage[font=small,labelfont=bf]{caption}

\usepackage[colorlinks=true,linkcolor=blue,citecolor=blue,urlcolor=blue]{hyperref}

\numberwithin{equation}{section}

\theoremstyle{plain}
\newtheorem{theorem}{Theorem}[section]

\newtheorem{lemma}[theorem]{Lemma}
\newtheorem{corollary}[theorem]{Corollary}

\newtheorem{definition}[theorem]{Definition}
\newcommand{\Pro}{\par\noindent\rm\textsc{Proof. }}

\def\bX{{\boldsymbol{X}}}
\def\bY{{\boldsymbol{Y}}}

\def\dd{{{\hskip1pt}\rm{d}}}

\def\V{{\mathcal V}}
\def\DC{{\mathcal R}}

\def\E{\mathbb{E}}
\def\C{\mathbb{C}}
\def\N{\mathbb{N}}
\def\R{\mathbb{R}}

\def\BState{\State\hskip-\ALG@thistlm}

\begin{document}

\title{\Large\textbf{A Basic Treatment of the Distance Covariance}}

\author{{Dominic Edelmann}\thanks{
Division of Biostatistics, German Cancer Research Center, Im Neuenheimer Feld 280, 69120 Heidelberg, Germany.  E-mail address: \href{mailto:dominic.edelmann@dkfz-heidelberg.de}{dominic.edelmann@dkfz-heidelberg.de}.
\endgraf
\ $^*$Division of Biostatistics, German Cancer Research Center, Im Neuenheimer Feld 280, 69120 Heidelberg, Germany.  E-mail address: \href{mailto:t.terzer@dkfz-heidelberg.de}{t.terzer@dkfz-heidelberg.de}.
\endgraf
\ $^\dag$Department of Statistics, Pennsylvania State University, University Park, PA 16802, U.S.A. E-mail address: \href{mailto:richards@stat.psu.edu}{richards@stat.psu.edu}.
\endgraf
\endgraf
\endgraf
} , {Tobias Terzer}$^*$, 
 and {Donald Richards}$^\dag$
 \endgraf
}

\maketitle

\begin{abstract}
The distance covariance of Sz\'ekely, \textit{et al.} \cite{szekely2007} and Sz\'ekely and Rizzo \cite{szekely2009brownian}, a powerful measure of dependence between sets of multivariate random variables, has the crucial feature that it equals zero if and only if the sets are mutually independent.  Hence the distance covariance can be applied to multivariate data to detect arbitrary types of non-linear associations between sets of variables.  We provide in this article a basic, albeit rigorous, introductory treatment of the distance covariance.  Our investigations yield an approach that can be used as the foundation for presentation of this important and timely topic even in advanced undergraduate- or junior graduate-level courses on mathematical statistics.  
\end{abstract}

\section{Introduction}
\label{sec_introduction}

\begin{figure}[!t]
\captionsetup{width=0.875\textwidth}
\includegraphics[width=\textwidth]{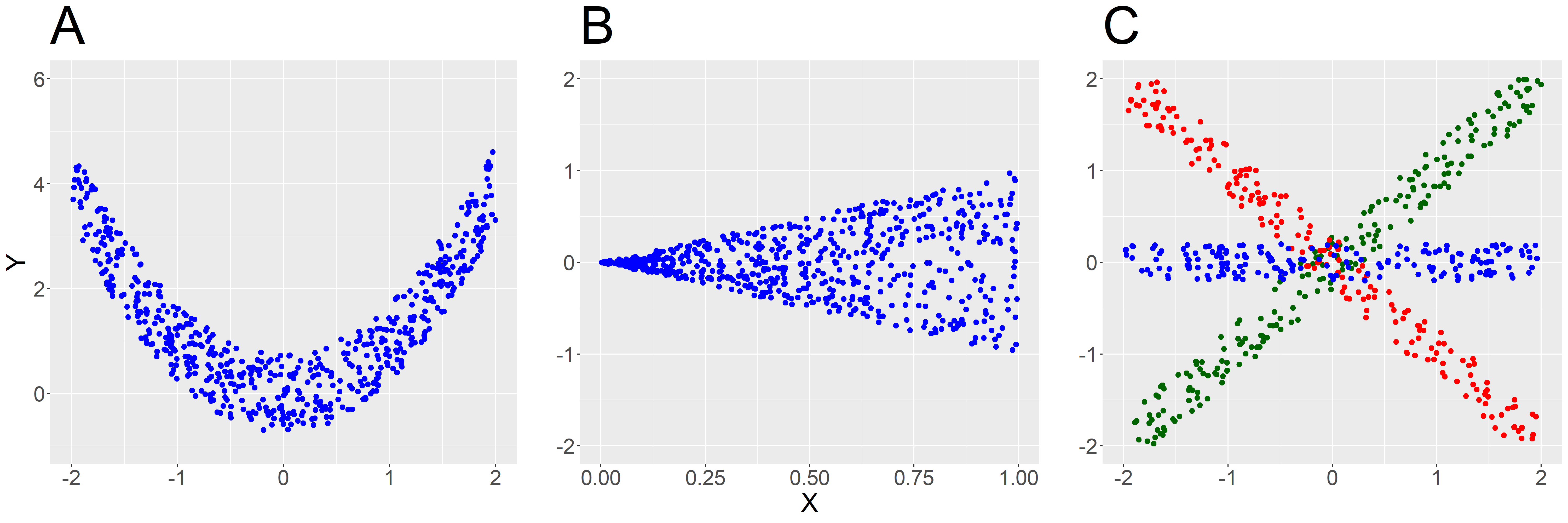}
\caption{The sub-figures {\bf A-C} represent scatter-plots of bivariate samples $(\bX,\bY)$ with $n=600$ data points to which independence tests, based on the distance covariance and classical covariance, were applied. In each case a distance covariance permutation test using $100,000$ permutations yields $p$-values of $10^{-5}$, demonstrating that the distance covariance is able to detect these dependencies. The $p$-values of permutation tests based on the classical covariance with $100,000$ permutations are $0.663$, $0.129$, and $0.889$ for {\bf A}, {\bf B}, and {\bf C}, respectively.}
\label{fig:demo}
\end{figure}

The distance covariance, a measure of dependence between multivariate random variables $X$ and $Y$, was introduced by Sz\'ekely, Rizzo, and Bakirov \cite{szekely2007} and has since received extensive attention in the statistical literature.  A crucial feature of the distance covariance is that it equals zero if and only if $X$ and $Y$ are mutually independent. Hence the distance covariance is sensitive to arbitrary dependencies; this is in contrast to the classical covariance, which is generally capable of detecting only linear dependencies. This property is illustrated in Figure \ref{fig:demo}, which illustrates that tests based on the distance covariance are able to detect numerous types of non-linear associations even when tests based on the classical covariance may fail to detect many such statistical relationships.

\begin{figure}[!t]
	\captionsetup{width=0.875\textwidth}
	\includegraphics[width=\textwidth]{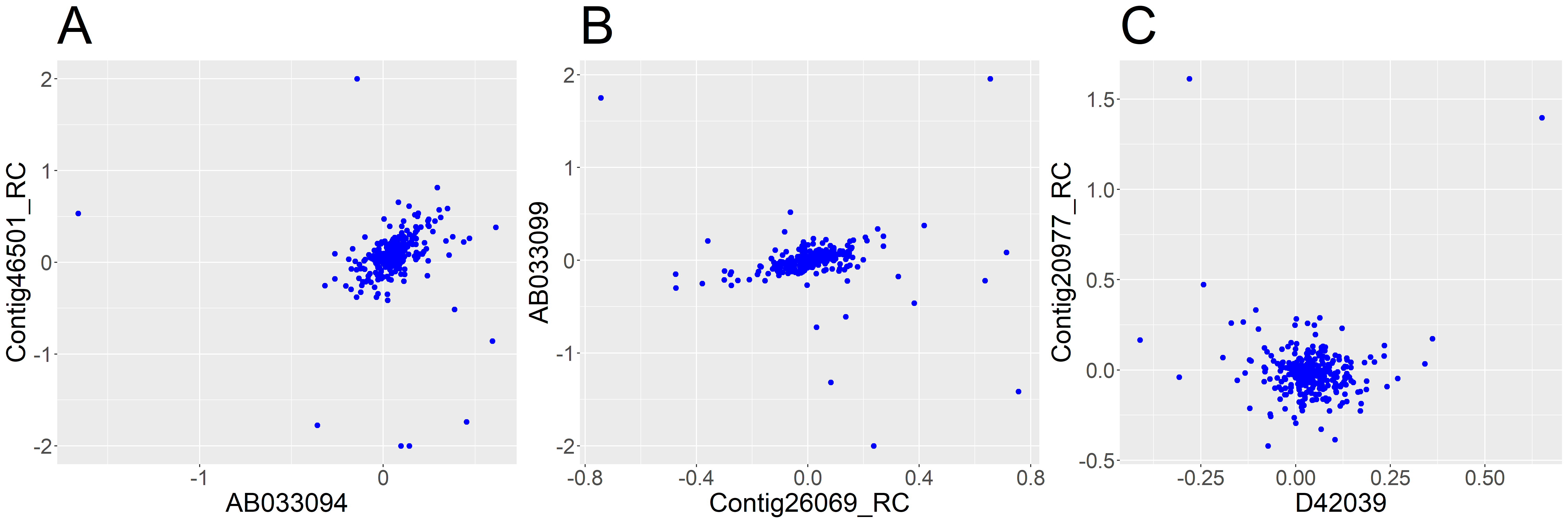}
	\caption{Sub-figures {\bf A-C} represent three scatter-plots of the expression values of genes in a breast cancer data set provided by Van De Vijver, \textit{et al.} \cite{van2002gene} ($n=295$ samples) on which permutation tests, based on the distance covariance and classical covariance, were applied. The $p$-values of the distance covariance permutation tests using $100,000$ permutations are {\bf A}: $10^{-5}$ ; {\bf B}: $10^{-5}$; {\bf C}: $3.00 \times 10^{-4}$  . The $p$-values of permutation tests based on the classical covariance with $100,000$ permutations are {\bf A}: $0.079$ ; {\bf B}: $0.503$; {\bf C}: $0.930$.}
	\label{fig:vijver}
\end{figure}

While the dependencies illustrated in Figure \ref{fig:demo} clearly represent purely illustrative examples, the sensitivity of the distance covariance to arbitrary dependencies can be very useful for applications. This is demonstrated in Figure \ref{fig:vijver}, where we show three dependencies between expression values genes in the breast cancer data set by Van De Vijver, \textit{et al.} \cite{van2002gene}; all these dependencies can be detected by the distance covariance but not by the classical covariance.

For comparisons of the distance covariance and classical covariance in applications to data, see the examples given by Sz\'ekely and Rizzo \cite[Section 5.2]{szekely2009brownian} and Dueck, \textit{et al.} \cite[Section~5]{dueck2014affinely}; for extensive numerical experiments and fast algorithms for computing the distance covariance, see Huo and Sz\'ekely \cite[Section 5]{huo2016fast}.  We also refer to Sejdinovic, \textit{et al.} \cite{Sejdinovic2013}, Dueck, \textit{et al.} \cite{dueck2014affinely}, Sz\'ekely and Rizzo \cite{szekely2009brownian,szekely2014partial}, Huo and Sz\'ekely \cite{huo2016fast}, and Edelmann, \textit{et al.} \cite{edelmann2020distance}, representing only a few of the many authors who have given further theoretical results on the distance covariance and distance correlation coefficients; and to Zhou \cite{zhou2012}, Fiedler \cite{fiedler2016}, and Edelmann, et al. \cite{edelmann2019review} as among the applications to time series analyses.  Many applications to data analysis of the distance correlation coefficient and the distance covariance are now available, including: Kong, \textit{et al.} \cite{kong2012using} on data in sociology, Mart\'inez-G\'omez, textit{et al.} \cite{martinez2014distance} and Richards, \textit{et al.} \cite{richards2014interpreting} on astrophysical databases and galaxy clusters, Dueck, \textit{et al.} \cite{dueck2014affinely} on time series analyses of wind vectors at electricity-generating facilities, Richards \cite{richards2017} on the relationship between the strength of gun control laws and firearm-related homicides, Zhang, \textit{et al.} \cite{zhang2018} for remote sensing applications, and Ohana-Levi, \textit{et al.} \cite{ohana2020} on grapevine transpiration.

The original papers of Sz\'ekely, \textit{et al.} \cite{szekely2007} and Sz\'ekely and Rizzo \cite{szekely2009brownian} are now widely recognized as seminal and important contributions to measuring dependence between sets of random variables; however, the exposition therein includes some ingenious arguments that may make the material challenging to readers not having an advanced background in mathematical statistics.  With the benefit of hindsight, we are able to provide in this article a simpler, albeit mathematically rigorous, introduction to the distance covariance that can be taught even in an undergraduate-level course covering the basic theory of U-statistics.  Other than standard U-statistics theory and some well-known properties of characteristic functions, the requirements for our treatment are a knowledge of multidimensional integration and trigonometric inequalities, as covered in a course on undergraduate-level advanced calculus.  Consequently, we hope that this treatment will prove to be beneficial to non-mathematical statisticians.  

Our presentation introduces the distance covariance as an important alternative to the classical covariance.  Moreover, the distance covariance constitutes a particularly interesting example of a U-statistic since it includes both the ``non-degenerate'' and ``first-order degenerate'' cases of the asymptotic distribution theory of U-statistics, these corresponding to the situations in which $X$ and $Y$ are dependent, leading to the non-degenerate case, or $X$ and $Y$ are independent, leading to the first-order degenerate case of the asymptotic theory.  

Throughout the exposition, $\|\cdot\|$ denotes the Euclidean norm and $\langle \cdot,\cdot \rangle$ the corresponding inner product.  Also, we denote by $|\cdot|$ the modulus in $\C$ or the absolute value in $\R$, and the imaginary unit is $i = \sqrt{-1}$.  


\section{The fundamental integral of distance covariance theory}

Following Sz\'ekely, \textit{et al.} \cite{szekely2007}, we first establish a closed-form expression for an integral that plays a central work in this article, leading to the equivalence of two crucial expressions for the distance covariance.  The first expression displays the distance covariance as an integrated distance between the joint characteristic function of $(X,Y)$ and the product of the marginal characteristic functions of $X$ and $Y$; we will deduce from this expression that the distance covariance equals zero if and only if $X$ and $Y$ are independent.  The second expression allows us to derive consistent distance covariance estimators that are expressible as polynomials in the distances between random samples.

Since the ability to characterize independence and the existence of easily computable estimators are arguably the most important properties of the distance covariance, we will refer to this integral as the {\it fundamental integral of distance covariance}.

\begin{lemma}
\label{lem:fund}
For $x \in \R^p$, 
\begin{equation}
\label{eq_fund_integral}
\int_{\R^p} \frac{1-\cos\,\langle t,x \rangle}{\|t\|^{p+1}} \dd t = \frac{\pi^{(p+1)/2}}{\Gamma\big((p+1)/2\big)} \, \|x\|.
\end{equation}
\end{lemma}

\Pro
Since \eqref{eq_fund_integral} is valid for $x = 0$, we need only treat the case in which $x \neq 0$.  

Denote by $I_p$ the integral in \eqref{eq_fund_integral}.  For $p=1$, replacing $t$ by $t/x$ yields 
\begin{equation}
\label{c1_integral}
I_1 = \int_{-\infty}^\infty \frac{1-\cos t x }{t^{2}} \dd t 
= |x| \, \int_{-\infty}^\infty \frac{1-\cos t}{t^2} {\dd t}.
\end{equation}
Denoting the latter integral in \eqref{c1_integral} by $c_1$, it follows by integration-by-parts that 
\begin{equation}
\label{c1}
c_1 = 2\, \int_0^\infty \frac{1- \cos t}{t^2} {\dd t} = 2\, \int_0^\infty \frac{\sin t}{t} {\dd t} = \pi,
\end{equation}
the last equality being classical in calculus (Spivak \cite[Chapter 19, Problem 43]{spivak1994calculus}).  

For general $p$, note that $I_p$ is invariant under orthogonal transformations $H$ of $x$:
\begin{align*}
\int_{\R^p} \frac{1-\cos \, \langle t,Hx \rangle}{\|t\|^{p+1}} \dd t 
&= \int_{\R^p} \frac{1-\cos \, \langle Ht, Hx \rangle}{\|H \, t\|^{p+1}} \dd t \\
&= \int_{\R^p} \frac{1-\cos \, \langle t, x \rangle}{\|t\|^{p+1}} \dd t,
\end{align*}
where the first equality follows from the transformation $t \mapsto Ht$, which leaves the Lebesgue measure $\dd t$ unchanged; and the second equality holds because the norm and the inner product are orthogonally invariant.  Therefore, in evaluating $I_p$ we may replace $x$ by $\|x\|(1,0,\ldots,0)$; letting $t=(t_1,\ldots,t_p)$, we obtain 
\begin{equation}
\label{cp_integral}
I_p = \int_{\R^p} \frac{1-\cos \, (t_1 \|x\|)}{\|t\|^{p+1}} \dd t = \|x\| \int_{\R^p} \frac{1-\cos t_1}{\|t\|^{p+1}} \dd t,
\end{equation}
the last equality obtained by replacing $t_j$ by $t_j/\|x\|$, $j=1,\ldots,p$.  

Denoting by $c_p$ the latter integral in \eqref{cp_integral}, we substitute in that integral $t_j = v_j$, $j=1,\ldots,p-1$, and $t_p = p^{-1/2} (v_1^2 + \cdots + v_{p-1}^2)^{1/2} v_p$.  As the Jacobian of this transformation is $p^{-1/2} (v_1^2 + \cdots + v_{p-1}^2)^{1/2}$, we obtain 
\begin{align}
\label{student}
c_p &= p^{-1/2} \int_{\R^{p-1}} \frac{1-\cos v_1}{(v_1^2 + \cdots + v_{p-1}^2)^{p/2}} \dd v_1 \cdots \dd v_{p-1} \cdot \int_{-\infty}^\infty \frac{\dd v_p}{(1 + p^{-1} v_p^2)^{(p+1)/2}} \nonumber \\
&= p^{-1/2} c_{p-1} \, \int_{-\infty}^\infty \frac{\dd v_p}{(1 + p^{-1} v_p^2)^{(p+1)/2}}.
\end{align}
As the remaining integral in \eqref{student} is the familiar normalizing constant of the Student's $t$-distribution on $p$ degrees-of-freedom, we obtain 
$$
c_p = \frac{\pi^{1/2} \Gamma(p/2)}{\Gamma\big((p+1)/2\big)} \, c_{p-1}.
$$
Starting with $c_1 = \pi$, we solve this recursive equation for $c_p$, obtaining \eqref{eq_fund_integral}.  
$\qed$

\section{Two representations for the distance covariance} \label{sec:alter}

We now introduce the representations of the distance covariance mentioned above. Following Sz\'ekely, \textit{et al.} \cite{szekely2007}, we define the distance covariance through its characteristic function representation. For jointly distributed random vectors $X \in \R^p$ and $Y \in \R^q$, let $\phi_{X,Y} (s,t) = \E e^{i \langle s , X \rangle+ i \langle t , Y \rangle }$ be the joint characteristic function of $(X,Y)$ and $\phi_{X} (s) = \phi_{X,Y} (s,0)$ and $\phi_{Y} (t) = \phi_{X,Y} (0,t)$ be the corresponding marginal characteristic functions.

\begin{definition}
The {\it distance covariance} $\V(X,Y)$ between $X$ and $Y$ is defined as the nonnegative square-root of
\begin{equation}
\label{eq:dcov1}
\V^2(X,Y) = \frac{1}{c_p c_q} \int_{\R^p}\int_{\R^q} \frac{|\phi_{X,Y}(s,t) - \phi_X(s) \phi_Y (t)|^2}{\|s\|^{p+1} \, \|t\|^{p+1}} \, \dd s \, \dd t,
\end{equation}
where $c_p$ is the normalizing constant in \eqref{eq_fund_integral}. 	
\end{definition}   

As the integrand in \eqref{eq:dcov1} is nonnegative, it follows that $\V^2(X,Y) \geq 0$.  Further, we will show in Corollary \ref{cor:existence} that $\V^2(X,Y) < \infty$ whenever $X$ and $Y$ have finite first moments.  

An advantage of the representation \eqref{eq:dcov1} is that it directly implies one of the most important properties of the distance covariance, viz., the characterization of independence.

\begin{theorem}
For all $X$ and $Y$, $\V^2(X,Y) = 0$ if and only if $X$ and $Y$ are independent.
\end{theorem}

\Pro
If $X$ and $Y$ are independent then $\phi_{X,Y} (s,t) = \phi_X(s) \phi_Y (t)$ for all $s$ and $t$; hence $\V^2(X,Y) = 0$.

Conversely, if $X$ and $Y$ are not independent then the functions $\phi_{X,Y}(s,t)$ and $\phi_X(s) \phi_Y(t)$ are not identical (Van der Vaart \cite[Lemma 2.15]{van2000asymptotic}).  Since characteristic functions are continuous then there exists an open set $\mathcal{A} \subseteq \R^p \times \R^q$ such that $|\phi_{X,Y}(s,t) - \phi_X(s) \phi_Y(t)|^2 > 0$ for all $(s,t) \in \mathcal{A}$.  Hence, by \eqref{eq:dcov1}, $\V^2(X,Y) > 0$.  
$\qed$

\medskip

For the purpose of deriving estimators for $\V^2(X,Y)$, we now apply Lemma \ref{lem:fund} to obtain a second representation of the distance covariance.

\begin{theorem}
Suppose that $(X_1,Y_1),\ldots,(X_4,Y_4)$ are independent, identically distributed (i.i.d.) copies of $(X,Y)$.  Then
\begin{equation} 
\label{eq:dcov2}
\V^2(X,Y ) = \E \Big[\|X_1 - X_2 \| \cdot \|Y_1 - Y_2 \|- 2\|X_1 - X_2 \| \cdot \|Y_1 - Y_3 \| + \|X_1 - X_2 \| \cdot \|Y_3 - Y_4 \| \Big].
\end{equation}
\end{theorem}

\Pro
First, we observe that the numerator in the integrand in \eqref{eq:dcov1} equals 
\begin{align*}
|\phi_{X,Y}&(s,t) - \phi_X(s) \phi_Y (t)|^2 \\
&= (\phi_{X,Y} (s,t) - \phi_X(s) \phi_Y (t)) \, \overline{(\phi_{X,Y} (s,t) - \phi_X(s) \phi_Y (t))} \\
&= \E\big[e^ {i \langle s , X_1-X_2 \rangle+ i \langle t , Y_1-Y_2 \rangle} - 2\,e^ {i \langle s , X_1-X_2 \rangle+ i \langle t , Y_1-Y_3 \rangle} + \,e^ {i \langle s , X_1-X_2 \rangle+ i \langle t , Y_3-Y_4 \rangle}\big].
\end{align*}
Since the latter expression is real, any term of the form $e^{i z}$, $z \in \R$, can be replaced by $\cos z$.  Hence, by \eqref{eq:dcov1}, 
\begin{align}
\label{eq:indecomp1}
c_p c_q \V^2(X,Y) =\int_{\R^p}\int_{\R^q} \frac{A_{12}(s,t)- 2\,A_{13}(s,t) + A_{34}(s,t)} {\|s\|^{p+1} \| t \| ^{q+1}} \dd s \dd t
\end{align}
where, for each $(j,k)$, 
\begin{equation}
\label{Ajk}
A_{jk}(s,t) = \E \cos{\big(\langle s, X_1 - X_2 \rangle +  \langle t ,Y_j - Y_k \rangle \big)}.
\end{equation}
Replacing $t$ by $-t$ in \eqref{eq:indecomp1}, we also obtain 
\begin{align}
\label{eq:indecomp2}
c_p c_q \V^2(X,Y) =\int_{\R^p}\int_{\R^q} \frac{A_{12}(s,-t)- 2\,A_{13}(s,-t) + A_{34}(s,-t)} {\|s\|^{p+1} \| t \| ^{q+1}} \dd s \dd t,
\end{align}
and by adding \eqref{eq:indecomp1} and \eqref{eq:indecomp2}, we find that 
$$
c_p c_q \V^2(X,Y) = \int_{\R^p}\int_{\R^q} \frac{B_{12}(s,t)- 2\,B_{13}(s,t) + B_{34}(s,t)} {\|s\|^{p+1} \| t \| ^{q+1}}  \dd s \dd t 
$$
where for each $(j,k)$, 
\begin{equation}
\label{Bjk}
B_{jk}(s,t) = \frac{1}{2} \big( A_{jk}(s,t) + A_{jk}(s,-t) \big).
\end{equation}
On applying to \eqref{Ajk} and \eqref{Bjk} the trigonometric identity, 
$$
\cos(x+y) + \cos(x-y) = 2\, \cos x \cos y,
$$
we deduce that 
\begin{equation}
\label{Bjk2}
B_{jk}(s,t) = \E \big[ \cos{\langle s, X_1 - X_2 \rangle} \,  \cos{\langle t, Y_j - Y_k \rangle} \big].
\end{equation}

For $j,k, \in \{1,2,3,4\}$, we apply to \eqref{Bjk2} the elementary identity, 
\begin{align}
\label{eq:oneminus}
\cos{\langle s, X_1 - X_2 \rangle} \,  \cos{\langle t, Y_j - Y_k \rangle} &=
\big(1-\cos{\langle s, X_1 - X_2 \rangle}\big) \, \big(1-\cos{\langle t, Y_j - Y_k \rangle}\big) \nonumber \\
& \quad\quad - 1 +  \cos{\langle s, X_1 - X_2 \rangle} + \cos{\langle t, Y_j - Y_k \rangle};
\end{align}
then we obtain 
\begin{align*}
c_p & c_q \V^2(X,Y) \\
&= \int_{\R^p}\int_{\R^q} \Big(\E \big[\big(1-\cos{\langle s, X_1 - X_2 \rangle} \big) \,  \big(1-\cos{\langle t, Y_1 - Y_2 \rangle} \big) \big] \\
& \qquad\qquad\quad - 2 \, \E \big[\big(1-\cos{\langle s, X_1 - X_2 \rangle} \big) \, \big(1-\cos{\langle t, Y_1 - Y_3 \rangle} \big)\big] \\
& \qquad\qquad\quad + \E \big[\big(1-\cos{\langle s, X_1 - X_2 \rangle} \big) \, \big(1-\cos{\langle t, Y_3 - Y_4 \rangle} \big)\big] \Big) \frac{\dd s \dd t}{\|s\|^{p+1} \|t\|^{q+1}},
\end{align*}
which is obtained by decomposing all summands on the right-hand side using Eq. \eqref{eq:oneminus} and observing that all terms which are not of the form $ \E[ \cos{\langle s, X_i - X_j \rangle} \,  \cos{\langle t, Y_l - Y_k \rangle}]$ cancel each other.  By applying the Fubini-Tonelli Theorem and the linearity of expectation and integration, we obtain 
\begin{align*}
c_p & c_q \V^2(X,Y) \\
&= \E \int_{\R^p}\int_{\R^q} \Big[ \big(1-\cos{\langle s, X_1 - X_2 \rangle} \big) \, \big(1-\cos{\langle t, Y_1 - Y_2 \rangle} \big) \\
& \qquad\qquad\qquad - 2 \, \big(1-\cos{\langle s, X_1 - X_2 \rangle} \big) \, \big(1-\cos{\langle t, Y_1 - Y_3 \rangle} \big) \\
& \qquad\qquad\qquad + \big(1-\cos{\langle s, X_1 - X_2 \rangle} \big) \, \big(1-\cos{\langle t, Y_3 - Y_4 \rangle} \big) \Big] {\hskip-1.54312pt} \frac{\dd s \dd t}{\|s\|^{p+1} \| t \| ^{q+1}}.
\end{align*}
The proof is completed by applying Lemma \ref{lem:fund} to calculate these three integrals.  
$\qed$

\medskip

Before establishing estimators for $\V^2(X,Y)$, we remark briefly on the assumptions necessary for the existence of the distance covariance.

\begin{corollary} {\label{cor:existence}}
	Suppose that $\E \|X\| < \infty$ and $\E \|Y\| < \infty$.  Then $\V^2(X,Y) < \infty$.
\end{corollary}

\Pro
From the representation \eqref{eq:dcov2}, we directly obtain the alternative representation
	\begin{multline}
\V^2(X,Y ) = \E \Big[\|X_1 - X_2 \| \, \|Y_1 - Y_2 \|- \|X_1 - X_2 \| \, \|Y_1 - Y_3 \| \\ - \|X_1 - X_2 \| \, \|Y_2 - Y_3 \|   + \|X_1 - X_2 \| \, \|Y_3 - Y_4 \| \Big].
	\end{multline}
Applying the triangle inequality yields 
$$
\|X_1 - X_2 \| \, \|Y_1 - Y_2 \|- \|X_1 - X_2 \| \, \|Y_1 - Y_3 \|  - \|X_1 - X_2 \| \, \|Y_2 - Y_3 \| \leq 0,
$$
and hence 
\begin{align*}
0 \leq \V^2(X,Y) &\leq \E \|X_1 - X_2 \| \, \|Y_3 - Y_4 \| \\
&= \E \|X_1 - X_2 \| \E \|Y_3 - Y_4 \| \leq 4 \, E \|X\| \E \|Y\|,
\end{align*}
where the last inequality follows again by the triangle inequality.  
$\qed$

\section{Asymptotic theory for estimating the distance covariance}

Using the representation of the distance covariance given in \eqref{eq:dcov2}, it is straightforward to derive a U-statistic estimator for $\V^2(X)$.  Specifically, we define the symmetric kernel function
\begin{multline}
\label{eq:kernel}
h\big((X_1,Y_1),\ldots,(X_4,Y_4)\big) \\ = \frac{1}{24} \sum
\big( \|X_i - X_j \| \,  \|Y_i - Y_j \|  - 2 \, \|X_i - X_j \| \,  \|Y_i - Y_k \| + \|X_i - X_j \| \,  \|Y_k - Y_l \| \big),
\end{multline}
where the sum is over all $i,j,k,l \in \{1,2,3,4 \}$ such that $i$, $j$, $k$, and $l$ are distinct.  

It follows from the representation \eqref{eq:dcov2} that each of the $24$ summands in \eqref{eq:kernel} has expectation $\V^2(X,Y)$.  Therefore, 
	$$
		\E h\big((X_1,Y_1),\ldots,(X_4,Y_4)\big)  = \V^2(X,Y).
	$$
Letting $(X_1,Y_1),\ldots,(X_n,Y_n)$ be a random sample from $(X,Y)$, we find that an unbiased estimator of $\V^2(X,Y)$ is 
\begin{equation} \label{eq:est}
\widehat{\Omega} = {\binom{n}{4}}^{-1} \sum_{1 \le i < j < k < l \le n} h\big((X_i,Y_i),(X_j,Y_j),(X_k,Y_k),(X_l,Y_l)\big).
\end{equation}
We can now derive the consistency and asymptotic distribution of this estimator using standard U-statistic theory (Lee \cite{lee2019u}). For this purpose, let us define
$$
h_1(x,y)  = \E\big[h\big((x,y),(X_2,Y_2),(X_3,Y_3),(X_4, Y_4)\big)\big].
$$
and
$$
h_2((x_1,y_1),(x_2,y_2))  = \E\big[h\big((x_1,y_1),(x_2,y_2),(X_3,Y_3),(X_4, Y_4)\big)\big].
$$
The preceding formulas and a classical result on U-statistics (Hoeffding \cite[Theorem 7.1]{Hoeffding1948}) leads immediately to a proof of the following result.  

\begin{theorem} 
\label{th:asympnormal}
Suppose that $0 < \mbox{Var}\,(h_1(X,Y)) < \infty$.  Then $\sqrt{n} \big(\widehat{\Omega} - \V^2(X,Y)\big) \stackrel{P}{\longrightarrow} Z$ as $n \to \infty$, where $Z \sim \mathcal{N}\big(0,16 \mbox{Var}\,(h_1(X,Y)\big)$.
\end{theorem}

\medskip

Except for pathological examples, Theorem \ref{th:asympnormal} provides the asymptotic distribution of $\V^2(X,Y)$ if $X$ and $Y$ are dependent.  For the crucial case of independent $X$ and $Y$, however, the asymptotic distribution of $\sqrt{n} (\widehat{\Omega} - \V(X,Y)^2)$ is degenerate; in this case, the asymptotic distribution can be derived using results on first-order degenerate U-statistics (Lee \cite[Section 3.2.2]{lee2019u}).

\begin{lemma} 
\label{lem:h1h2}
Let $X$ and $Y$ be independent, and $(X_1,Y_1)$ and $(X_2,Y_2)$ be i.i.d.~copies of $(X,Y)$.  Then
$h_1(x,y) \equiv 0$ and $\mbox{Var}\,\big(h_2((X_1,Y_1),(X_2,Y_2))\big) = \V^2(X,X) \, \V^2(Y,Y)/36$.
\end{lemma}

\medskip

The proof follows by elementary, but lengthy, transformations and may be left as an exercise to students.  A complete proof is provided by Huang and Ho \cite[Appendices B.6 and B.7]{huang2017statistically}.  

Finally, the following result follows directly from Lemma \ref{lem:h1h2} and classical results on the distributions of first-order degenerate U-statistics (Lee \cite[Section 3.2.2]{lee2019u}).  

\begin{theorem}
Let $X$ and $Y$ be independent, with $\E(\|X\|) < \infty$ and $\E(\|Y\|) < \infty$.  Then, 
\begin{align}
n \, \big(\widehat{\Omega} - \V^2(X,Y)\big) \stackrel{\mathcal{D}}{\longrightarrow} 6 \, \sum_{i=1}^\infty \lambda_i (Z_i^2 - 1),
\end{align}
as $n \to \infty$, where $Z_1,Z_2,\ldots$ are i.i.d.~standard normal random variables and $\lambda_1,\lambda_2,\ldots$ are the eigenvalues of the integral equation
$$
\E \big[h_2\big((x_1,y_1),(X_2,Y_2)\big) \, f(X_2,Y_2)\big] = \lambda f(x_1,y_1).
$$	
\end{theorem}

\section{Concluding Remarks}

In this article, we have derived under minimal technical requirements the most important statistical properties of the distance covariance.  From this starting point, there are several additional interesting topics that can be explored, e.g., as instructional assignments:

\medskip

\noindent
(i) The estimator \eqref{eq:est} is $O(n^4)$ and is computationally inefficient. A straightforward combinatorial computation shows that an $O(n^2)$ estimator of $\V$ is given by 
\begin{equation}
\begin{aligned}
\widetilde{\Omega} = \frac{1}{n \, (n-3)}\Bigg[ & \sum_{i,j=1}^n \|X_i-X_j\| \|Y_i-Y_j\| \\
& \ +  \frac{1}{(n-1) \, (n-2)} \sum_{i,j=1}^n \|X_i-X_j\| \cdot \sum_{i,j=1}^n \|Y_i-Y_j\| \\ 
& \ - \frac{2}{(n-2)} \sum_{i,j,k=1}^n \|X_i-X_j\| \|Y_i-Y_k\|\Bigg];
\label{eq:ugendcov}
\end{aligned}
\end{equation}
see Huo and Sz\'ekely \cite{huo2016fast}.  

\medskip

\noindent
(ii) We remark that although no assumption was provided in Theorem \ref{th:asympnormal} to ensure that $\mbox{Var}(h_1(X,Y)) < \infty$, it can be shown that this condition holds whenever $X$ and $Y$ have finite second moments; see Edelmann, \textit{et al.} \cite{edelmann2020distance}.

\medskip

\noindent
(iii) Important contributions of Sz\'ekely, \textit{et al.} \cite{szekely2007} and Sz\'ekely and Rizzo \cite{szekely2009brownian} are based on the {\it distance correlation coefficient}, which is defined as the nonnegative square-root of
$$
\DC^2(X,Y) = \frac{\V^2(X,Y)}{\sqrt{\V^2(X,X) \V^2(Y,Y)}}.
$$
Numerous properties of $\DC^2(X,Y)$ (see, e.g., Sz\'ekely, \textit{et al.} \cite[Theorem 3]{szekely2007}) may be derived using the methods that we have presented here.	

\medskip

We also remark on the fundamental integral, \eqref{eq_fund_integral}, that underpins the entire distance covariance and distance correlation theory.  As noted by Dueck, \textit{et al.} \cite{dueck2015}, the fundamental integral and variants of it have appeared in functional analysis (Gelfand and Shilov \cite[pp.~192--195]{GelfandShilov64}), in Fourier analysis (Stein \cite[pp.~140~and~263]{Stein70}), and in the theory of fractional Brownian motion on generalized random fields (Chil`es and Delfiner \cite[p.~266]{ChilesDelfiner12}; Reed, \textit{et al.} \cite{ReedLeeTruong95}).  

The fundamental integral also extends further.  For $m \in \N$ and $v \in \R$, define 
\begin{equation}
\label{truncated_cosine}
\cos_m(v) := \sum_{j=0}^{m-1} (-1)^j \frac{v^{2j}}{(2j)!},
\end{equation}
the truncated Maclaurin expansion of the cosine function.  Dueck, \textit{et al.} \cite{dueck2015} proved that for $\alpha \in \C$, 
\begin{equation}
\label{eq_fund_integral_2}
\int_{\R^d}\frac{\cos_m(\langle t,x\rangle) - \cos(\langle t,x\rangle)}{\|t\|^{d+\alpha}} \,\dd t = 
\frac{2\pi^{p/2} \, \Gamma(1-\alpha/2)}{\alpha \, 2^{\alpha} \, \Gamma\big((p+\alpha)/2\big)} \, \|x\|^{\alpha},
\end{equation}
with absolute convergence if and only if $2(m-1) < \Re(\alpha) < 2m$.  For $m=1$ and $\alpha = 1$, \eqref{eq_fund_integral_2} reduces to \eqref{eq_fund_integral}.  Further, for $m=1$ and $0 < \alpha < 2$, the integral \eqref{eq_fund_integral_2} provides the L{\'e}vy-Khintchine representation of the negative definite function $\|x\|^\alpha$, thereby linking the fundamental integral to the probability theory of the stable distributions.  

In conclusion, the statistical analysis of data through distance covariance and distance correlation theory, by means of the fundamental integral, is seen to be linked closely to many areas of the mathematical sciences.  

\bigskip
\bigskip

\noindent
{\bf Acknowledgements}.  D. Edelmann gratefully acknowledges financial support from the Deutsche Forschungsgemeinschaft (Grant number 417754611).  The authors are grateful to the editor of the special issue and to a reviewer for comments on the manuscript.

\end{document}